\documentclass[10pt,twoside]{article}
\usepackage[paperwidth=510.236pt,paperheight=737.008pt,top=55pt,bottom=55pt,left=52pt,right=52pt]{geometry}
\usepackage[T1]{fontenc}
\usepackage{amsmath,amssymb,amsthm,mathtools}
\usepackage{newtxtext,newtxmath}
\usepackage{graphicx}
\usepackage{enumitem}
\usepackage{fancyhdr}
\usepackage{microtype}
\usepackage{caption}
\usepackage{float}
\usepackage{needspace}
\usepackage[hidelinks]{hyperref}
\usepackage{xcolor}

\newtheoremstyle{acta}
{5pt}{5pt}
{\itshape}
{}
{\bfseries}
{}
{0.65em}
{}
\theoremstyle{acta}
\newtheorem{theorem}{Theorem}[section]
\newtheorem{lemma}[theorem]{Lemma}
\newtheorem{proposition}[theorem]{Proposition}
\newtheorem{corollary}[theorem]{Corollary}
\newtheorem{remark}[theorem]{Remark}
\newtheorem{example}[theorem]{Example}

\newcommand{\CP}{\mathbb{CP}}
\newcommand{\RP}{\mathbb{RP}}
\newcommand{\Z}{\mathbb{Z}}
\newcommand{\R}{\mathbb{R}}
\newcommand{\C}{\mathbb{C}}
\newcommand{\id}{\mathrm{id}}
\newcommand{\Pin}{\operatorname{Pin}}
\newcommand{\Diff}{\operatorname{Diff}}
\newcommand{\diag}{\operatorname{diag}}
\newcommand{\xt}{\operatorname{xt}}
\newcommand{\Int}{\operatorname{int}}

\begin{document}
\thispagestyle{empty}

\noindent\begin{minipage}[t]{0.58\textwidth}
\small Acta Mathematica Sinica, English Series\\
Revised manuscript: August 21, 2026
\end{minipage}\hfill
\begin{minipage}[t]{0.37\textwidth}
\raggedleft\small Acta Mathematica Sinica,\\English Series
\end{minipage}

\vspace{18pt}
\begin{center}
{\Large\bfseries On Closed Six-Manifolds Admitting Riemannian Metrics with\\[2pt]
Positive Sectional Curvature and Non-Abelian Symmetry\par}
\vspace{15pt}
{\large Yu Hang LIU\par}
\vspace{7pt}
{\small Department of Foundational Mathematics, School of Mathematics and Physics,\\
Xi'an Jiaotong-Liverpool University, Suzhou 215123, P. R. China\\
E-mail: Yuhang.Liu02@xjtlu.edu.cn\par}
\end{center}

\vspace{13pt}
\noindent\textbf{Abstract}\quad
We study the topology of closed, simply-connected, 6-dimensional Riemannian manifolds of positive sectional curvature which admit isometric actions by $SU(2)$ or $SO(3)$. We show that their Euler characteristic agrees with that of the known examples, i.e., $S^6$, $\CP^3$, the Wallach space $SU(3)/T^2$ and the biquotient $SU(3)//T^2$. We determine the non-trivial principal isotropy actions, classify the fixed-point-free $SU(2)$-actions without exceptional orbits, and give restrictions on the exceptional strata. In particular, an exceptional isotropy group for an $SU(2)$-action is cyclic of odd order, while tetrahedral exceptional isotropy occurs for the irreducible linear $SO(3)$-action on the round $6$-sphere.


\smallskip
\noindent\textbf{MR(2010) Subject Classification}\quad 53C20, 53C21, 53C23

\section{Introduction}

The study of Riemannian manifolds with positive sectional curvature is an old and fundamental subject in Riemannian geometry. There are very few compact examples of positively curved manifolds besides the so-called Compact Rank One Symmetric Spaces, which we will abbreviate as CROSS. In fact, the only further known examples occur only in dimension less than or equal to 24 and consist of homogeneous spaces \cite{BB76,Wal72}, biquotients \cite{Baz96,Esc82}, and one cohomogeneity one manifold in dimension 7 \cite{Dea11,GVZ11}.

The fundamental group of a compact Riemannian manifold with positive sectional curvature is finite, and it is trivial or equal to $\Z/2\Z$ in even dimensions. Furthermore, odd-dimensional positively curved closed manifolds are orientable. However, for simply connected closed manifolds, no general topological obstructions are known to separate the class of positively curved manifolds from the class of non-negatively curved manifolds, although there are many examples known to admit non-negative curvature.

There are several classification results for positively curved manifolds in low dimensions, though all of which require some ``symmetry'' conditions on the metric. Positively curved 3-manifolds are space forms \cite{Ham82}. In dimension 4, Hsiang and Kleiner showed that positively curved simply connected 4-manifolds with $S^1$ symmetry are homeomorphic to the 4-sphere $S^4$ or the projective space $\CP^2$ in \cite{HK89}; later Grove and Wilking improved the result to equivariant diffeomorphism in \cite{GW14}. In dimension 5, Xiaochun Rong showed that a $T^2$-invariant simply connected closed 5-manifold is homeomorphic to a 5-sphere in \cite{Ron02}. By the work of Barden \cite{Bar65} and Smale \cite{Sma62}, we know that there are no exotic 5-spheres, and thus such a 5-manifold is actually diffeomorphic to the standard 5-sphere.

Inspired by Hsiang and Kleiner's work, Karsten Grove proposed what is now called the ``symmetry program'' in \cite{Gro02}, which is to study positively curved manifolds with ``large'' symmetry group. Here ``large'' can have several different meanings. Many results were obtained in this direction, particularly for torus actions. For example, Grove and Searle proved the Maximal Rank theorem in \cite{GS94}, which states that the symmetry rank of an $n$-dimensional positively curved closed manifold is at most $[(n+1)/2]$, and in the case of equality the manifold is diffeomorphic to a sphere, $\RP^n$, $\CP^n$ or a lens space. In \cite{FR05}, Fuquan Fang and Xiaochun Rong showed that a closed simply connected $2m$-manifold $(m\geq 5)$ of positive sectional curvature on which an $(m-1)$-torus acts isometrically is homeomorphic to a complex projective space if and only if its Euler characteristic is not 2. Burkhard Wilking showed that when $n\geq 10$ and a positively curved closed simply connected $n$-manifold $M$ has symmetry rank at least $n/4+1$, $M$ is homeomorphic to $S^n$ or $\mathbb{HP}^{n/4}$ or homotopy equivalent to $\CP^{n/2}$ in \cite{Wil03}. Recently Kennard, Wiemeler and Wilking claimed that an even-dimensional positively curved manifold with $T^5$-symmetry has Euler characteristic at least 2 \cite{Ken20}.

After these classification results of positively curved manifolds invariant under torus actions, it is natural to investigate metrics with non-abelian symmetry. Wilking studied positively curved manifolds with high symmetry degree or low cohomogeneity relative to the dimension in \cite{Wil06}. Since all non-abelian compact Lie groups contain a rank 1 subgroup, $SU(2)$ or $SO(3)$, it is natural to study metrics invariant under $SU(2)$ or $SO(3)$. In dimension 5, Fabio Simas obtained a partial classification of positively curved 5-manifolds invariant under $SU(2)$ or $SO(3)$ in \cite{Sim16}. We point out here that Fabio Simas listed $SU(3)/SO(3)$ with the linear $SU(2)$-action as a candidate for positive curvature. But we note that this is not possible, since the fixed point set $(SU(3)/SO(3))^{\Z/2\Z}$ is diffeomorphic to $U(2)/O(2)$, which does not admit positive curvature.

In this paper we study 6-dimensional positively curved manifolds with $SU(2)$ or $SO(3)$ symmetry. This is also the first dimension where new examples other than the CROSSes have been constructed, which need to be recognized. They are the Wallach space $SU(3)/T^2$, where $T^2$ is the maximal torus, and the biquotient $SU(3)//T^2$, where
\[
T^2=\{(\diag(z,w,zw),\diag(1,1,\bar z^2\bar w^2))\mid z,w\in S^1\}\subset S(U(3)\times U(3))
\]
acts freely on $SU(3)$. On the first space one has an action by both $SO(3)$ and $SU(2)$, isometric in the positively curved metric, and on the second space an action by $SU(2)$ which commutes with $\diag(1,1,\bar z^2\bar w^2)$.

Our first result is:

\begin{theorem}\label{thm:main}
Let $M=M^6$ be a 6-dimensional closed simply connected Riemannian manifold of positive sectional curvature such that $SU(2)$ or $SO(3)$ acts isometrically and effectively on $M$. Then:
\begin{enumerate}[label=(\arabic*),leftmargin=2.2em]
\item The Euler characteristic $\chi(M)=2,4,6$.
\item The principal isotropy group is trivial unless $G=SO(3)$ and $M$ is equivariantly diffeomorphic to one of the two linear actions on $S^6$ given by $\R^3\oplus\R^4_{\mathrm{triv}}$ or $V_5\oplus\R^2_{\mathrm{triv}}$, where $V_5$ is the irreducible 5-dimensional real representation of $SO(3)$.
\item When the principal isotropy is trivial, every exceptional stratum has dimension at most one, and the closure of every connected component of the exceptional set contains a singular stratum. If $G=SU(2)$, every exceptional isotropy group is cyclic of odd order. If $G=SO(3)$, an exceptional isotropy group is cyclic, dihedral, tetrahedral, octahedral or icosahedral. The tetrahedral group occurs for a linear action on the round $S^6$.
\item Suppose $M/G=S^3$. Let $a$ denote the number of singular orbits with connected isotropy and let $b$ denote the number of singular orbits with disconnected isotropy. Then
\[
a+b\leq 3,\qquad \chi(M)=2a+b.
\]
In particular, $b$ is even. For $G=SU(2)$ one has $b=0$.
\end{enumerate}
\end{theorem}

Notice that in the known examples, one has indeed $\chi(M)=2,4,6$. Before stating the next theorems, we mention that the orbit space $M/G$ is homeomorphic to a 3-sphere or a 3-ball (see Theorem~\ref{thm:orbitspace}) unless $M$ is equivariantly diffeomorphic to $S^6$ with the fixed point homogeneous linear $SO(3)$-action in Example~\ref{ex:s6low}(2).

In the case of $G=SU(2)$ we will show:

\begin{theorem}\label{thm:su2main}
Let $M$ be a positively curved simply connected 6-manifold. Assume that $G=SU(2)$ acts on $M$ isometrically and effectively.
\begin{enumerate}[label=(\arabic*),leftmargin=2.2em]
\item If the fixed point set $M^G$ is non-empty, then $M$ is equivariantly diffeomorphic to a linear action on $S^6$ or $\CP^3$.
\item If $M^G$ is empty and the action has no exceptional orbits, then the action is equivariantly diffeomorphic to one of the actions on $S^6$, $S^2\times S^4$ or $SU(3)/T^2$ described in Remark~\ref{rem:su2models}.
\end{enumerate}
\end{theorem}

In the case of $G=SO(3)$ we have:

\begin{theorem}\label{thm:so3main}
Let $M$ be a positively curved simply connected 6-manifold. Assume that $G=SO(3)$ acts on $M$ isometrically and effectively and that the orbit space $M/G$ is a 3-ball whose boundary contains more than one orbit type, and that there are no exceptional orbits or interior singular orbits. Then $M^6$ is equivariantly homeomorphic to a linear action on $S^6$.
\end{theorem}

See Theorem~\ref{thm:so3ball} for further results in this special case. Explicit actions as in the above theorems are described in Section~\ref{sec:examples}.

The strategy to obtain these results is to analyze the structure of the orbit space and recover $M$ from $M/G$. We will show that $M/G$ is homeomorphic to $B^4$, $B^3$ or $S^3$. When $M/G=B^4$, the action is fixed point homogeneous and is completely understood. When $M/G=B^3$, we characterize the boundary and interior stratification for both $G=SU(2)$ and $G=SO(3)$. When $M/G=S^3$, we show that there are at most three singular orbits and give a formula for the Euler characteristic in terms of their isotropy groups. In all three cases, we describe the structure of singular orbit strata, which allows us to glue different pieces of singular orbits to recover the topology of $M$ if exceptional orbits do not occur. If exceptional orbits occur, we show that the closure of every exceptional component meets the singular set. The irreducible action of $SO(3)$ on the round $S^6$ also shows that polyhedral exceptional isotropy cannot be excluded.

This paper is organized as follows: Section 2 is devoted to preliminary knowledge on group actions and manifolds with positive sectional curvature. In Section 3 we study different orbit types for all possible actions considered in this paper. In Section 4 we prove Theorems~\ref{thm:main}, \ref{thm:su2main} and \ref{thm:so3main}, and study actions with more complicated orbit stratification. In Section~\ref{sec:examples} we describe the relevant isometric $SU(2)$- and $SO(3)$-actions on positively curved 6-manifolds.

\section{Preliminaries}

We start by recalling some basic definitions for group actions, see e.g. \cite{Bre72,AB15} for a reference. Let $G$ be a compact Lie group and $M$ be a compact smooth manifold. For a smooth action $\pi:G\times M\to M$, the $G$-orbit $G\cdot p$ through a point $p\in M$ is the submanifold $G\cdot p=\{gp\in M\mid g\in G\}$, the isotropy group or the stabilizer at $p\in M$ is defined as $G_p=\{g\in G\mid gp=p\}$, and we have $G\cdot p=G/G_p$. Furthermore, we denote the $G$-fixed point set by $M^G=\{p\in M\mid G\cdot p=p\}$. Note also that the fixed point set in an orbit has the form
\[
(G/K)^H=\{gK\in G/K\mid g^{-1}Hg\subset K\},
\]
where $H\subset K\subset G$. In particular, $(G/H)^H=N(H)/H$.

Points in the same $G$-orbits have conjugate isotropy groups. The isotropy type of a $G$-orbit $G/H$ is the conjugacy class of isotropy groups at points in $G/H$ and we denote it by $(H)$. We define $M_{(K)}$ to be the union of orbits with the same isotropy type $(K)$. For compact group actions on compact manifolds, there are only finitely many orbit types.

Among all orbit types of a given action, there exist maximal orbits $G/H$ with respect to inclusion of isotropy groups called the principal orbits. Non-principal orbits which have the same dimension as the principal orbit are called exceptional orbits, and orbits having lower dimension than principal ones are called singular orbits.

The orbit space $M^*=M/G$ is the union of its orbit strata $M^*_{(K)}=M_{(K)}/G$ which themselves are manifolds. The principal orbit stratum $M^*_{(H)}$ is an open, dense and connected subset of $M/G$. In particular the dimension of $M^*_{(H)}$ is called the cohomogeneity of the action. Codimension one strata in $M^*$ are called faces, which are part of $\partial M^*$. We will also use the fact that $M^K_{(K)}\to M^*_{(K)}$ is an $N(K)/K$-principal bundle and the structure group of $M_{(K)}\to M^*_{(K)}$ is $N(K)/K$.

The following theorem gives constraints on the exceptional orbits on simply-connected manifolds.

\begin{theorem}[{\cite[Theorem IV.3.12]{Bre72}}]\label{thm:nospecial}
Let $M$ be a simply-connected manifold and $G$ a compact group acting on $M$. Then $M^*$ is also simply connected and there are no exceptional orbits $G/K$ whose stratum $M^*_{(K)}$ has codimension one in $M^*$, the so-called special exceptional orbits.
\end{theorem}

The following theorem describes the manifold structure on $M^*$ for certain actions of cohomogeneity 3.

\begin{theorem}[{\cite[Corollary IV.4.7]{Bre72}}]\label{thm:3manifold}
Let $M$ be a compact and simply-connected manifold and $G$ a compact group acting on $M$. Suppose that all orbits are connected. Then $M^*$ is a simply connected topological 3-manifold with or without boundary.
\end{theorem}

For each orbit $G\cdot p$, let $T_p^\perp$ denote the normal space at $p$ to the orbit and $\nu_p$ the unit sphere in the normal space. $T_p^\perp$ admits a natural linear action by the isotropy group $G_p$, called the slice representation. The quotient $T_p^\perp/G_p$ is called the tangent cone of the orbit in the orbit space, and $\nu_p/G_p$ is the space of directions at $p$ and is denoted as $\Sigma_{[p]}$. If $G$ acts on $M$ by isometries, the orbit space, tangent cones and spaces of directions all inherit a metric from $M$. In particular, if we impose the positive curvature assumption on $M$, $M/G$ becomes an Alexandrov space with positive curvature.

We also note that $M_{(K)}\cap T_p^\perp=(T_p^\perp)^K$, and the slice theorem states that an equivariant neighborhood of $G/G_p$ has the form $G\times_{G_p}D(T_p^\perp)$, where $D(T_p^\perp)$ is a disk in $T_p^\perp$ and $G_p$ acts on $G$ via right multiplication and on $D(T_p^\perp)$ via the slice representation.

We now state the Extent Lemma. For any metric space $(X,d)$ and positive integer $q\geq 2$, we define the $q$-extent of $X$ as
\begin{equation}\label{eq:extent}
\xt_q(X)=\binom q2^{-1}\sup_{x_1,\ldots,x_q\in X}\sum_{1\leq i<j\leq q}d(x_i,x_j).
\end{equation}
In other words, $\xt_q(X)$ is the maximal average distance between points in $q$-tuples in $X$. When $q=2$, $\xt_2(X)$ is the diameter of $X$. The Extent Lemma from \cite{GS97} states that:

\begin{lemma}\label{lem:extent}
If $M/G$ is an Alexandrov space with positive curvature, then for all $(q+1)$-tuples $([x_0],\ldots,[x_q])$ in $M/G$, we have
\[
\frac1{q+1}\sum_{i=0}^q\xt_q(\Sigma_{[x_i]})>\frac\pi3.
\]
\end{lemma}

For the Euler characteristic we have:

\begin{theorem}\label{thm:eulerfixed}
\begin{enumerate}[label=(\arabic*),leftmargin=2.2em]
\item \cite{Kob58} If a torus $T$ acts smoothly on a closed smooth manifold $M$, then the Euler characteristic of $M$ equals that of $M^T$, that is, $\chi(M)=\chi(M^T)$.
\item \cite{PS02} If $M$ is a 6-dimensional simply connected Riemannian manifold with positive sectional curvature and $S^1$-symmetry, then $\chi(M)$ is positive and even.
\end{enumerate}
\end{theorem}

For totally geodesic submanifolds of positively curved manifolds, we have the so-called Connectedness Lemma due to Burkhard Wilking.

\begin{theorem}[Connectedness Lemma, {\cite{Wil03}}]\label{thm:connectedness}
Let $M^n$ be a compact $n$-dimensional Riemannian manifold with positive sectional curvature. Suppose that $N^{n-k}\subset M^n$ is a compact totally geodesic embedded submanifold of codimension $k$. Then the inclusion map $N^{n-k}\hookrightarrow M^n$ is $(n-2k+1)$-connected.
\end{theorem}

Recall that if we have a continuous map $f:X\to Y$ between two connected topological spaces $X$ and $Y$, and a positive integer $k$, then we say that $f$ is $k$-connected if $f_*$ induces isomorphisms on homotopy groups $\pi_i$, $1\leq i\leq k-1$, and a surjection on $\pi_k$.

For smooth actions on positively curved manifolds with nontrivial principal isotropy group, we have:

\begin{theorem}[Isotropy Lemma, {\cite{Wil06}}]\label{thm:isotropy}
Let $G$ be a compact Lie group acting isometrically and not transitively on a positively curved manifold $(M,g)$ with nontrivial principal isotropy group $H$. Then any nontrivial irreducible subrepresentation of the isotropy representation of $G/H$ is equivalent to a subrepresentation of the isotropy representation of $K/H$, where $K$ is an isotropy group such that the orbit stratum of $K$ is a boundary face in $M/G$ and $K/H$ is a sphere.
\end{theorem}

Let $G$ be a compact Lie group and $M$ a closed smooth manifold. For a smooth $G$-action on $M$ with non-empty fixed point set $M^G$, we define the fixed point cohomogeneity of this action as
\begin{equation}\label{eq:fixedcohom}
\operatorname{cohomfix}(M,G)=\dim(M/G)-\dim(M^G)-1,
\end{equation}
where $\dim(M^G)$ is the dimension of the fixed point component of largest dimension.

For Riemannian manifolds with positive sectional curvature and low fixed point cohomogeneity, we have the following classification.

\begin{theorem}[{\cite{GK04,GS97}}]\label{thm:fixedpointcohom}
If $M$ is a positively curved simply connected closed manifold which admits an isometric action by a compact group $G$ such that the fixed point cohomogeneity is less than or equal to 1, then $M$ is equivariantly diffeomorphic to a compact rank one symmetric space with an isometric $G$-action.
\end{theorem}

By examining the actions on rank one symmetric spaces in Section~\ref{sec:examples}, it follows that:

\begin{corollary}\label{cor:lowfixed}
Let $M$ be a 6-dimensional positively curved simply connected closed manifold which admits an isometric action by $G=SU(2)$ or $SO(3)$. If the action is fixed point homogeneous, then it is given by Examples~\ref{ex:s6low}(1), \ref{ex:s6low}(2), or \ref{ex:cp3}(1). If the action has fixed point cohomogeneity one, then it is given by Example~\ref{ex:s6low}(3).
\end{corollary}

We now state a version of the soul theorem in the setting of orbit spaces.

\begin{theorem}[Boundary Soul Lemma, {\cite[Theorem 1.2]{GK04}}]\label{thm:boundarysoul}
Let $M$ be a closed Riemannian manifold with positive sectional curvature and $G$ a compact Lie group acting isometrically on $M$. Suppose $M^*=M/G$ has nonempty boundary $\partial M^*$. Then:
\begin{enumerate}[label=(\arabic*),leftmargin=2.2em]
\item There exists a unique point $s_0\in M^*$, the soul of $M^*$, at maximal distance to $\partial M^*$.
\item The space of directions $\Sigma_{[s_0]}$ at $s_0$ is homeomorphic to $\partial M^*$.
\item The strata in $\Int(M^*)=M^*-\partial M^*$ belong to one of the following:
\begin{enumerate}[label=(\alph*),leftmargin=2.2em]
\item all of $\Int(M^*)$;
\item the soul point $s_0$;
\item a cone over strata in $\partial M^*$ with its cone point $s_0$ removed;
\item a stratum containing $s_0$ in its interior and whose boundary consists of strata in $\partial M^*$.
\end{enumerate}
\end{enumerate}
\end{theorem}

\begin{remark}\label{rem:soul}
We point out that the regular points in $M^*$ belong to either (3)(a) or (3)(c) in Theorem~\ref{thm:boundarysoul}. We also note that in \cite{GK04} it was claimed that the stratum in Theorem~\ref{thm:boundarysoul}(3)(d) is one-dimensional, but one easily gives examples where its dimension is higher. Consider the Hopf action by $S^1$ on $S^3$, and the $SO(3)$ action on $S^5$ given by the restriction of the diagonal action of $SO(3)$ on $S^5\subset\R^3\oplus\R^3$. Note that $S^5/SO(3)$ is a 2-dimensional hemisphere of radius $1/2$, which we denote by $H$. Now consider the product sum action by $S^1\times SO(3)$ on $S^3*S^n$, where $*$ is the spherical join. The orbit space is then $X=S^2(1/2)*H$, and the soul lies in a 2-dimensional stratum.
\end{remark}

We frequently use the knowledge of the subgroups of $SO(3)$ and $SU(2)$. For $SO(3)$ they are given by
\begin{itemize}[leftmargin=2em]
\item 0-dimensional subgroups: $\Z/k\Z$, $D_k$, $A_4$, $S_4$, $A_5$;
\item 1-dimensional subgroups: $SO(2)$, $O(2)$;
\end{itemize}
and for $SU(2)$ by
\begin{itemize}[leftmargin=2em]
\item 0-dimensional subgroups: $\Z/k\Z$, binary dihedral groups, inverse images of $A_4$, $S_4$, $A_5$ in $SU(2)$;
\item 1-dimensional subgroups: $U(1)$, $\Pin(2)=N(U(1))$.
\end{itemize}
Note that the only finite subgroups of $SU(2)$ which do not contain the center $\Z/2\Z$ are cyclic groups of odd order.

\begin{figure}[ht]
\centering
\includegraphics[width=.68\textwidth]{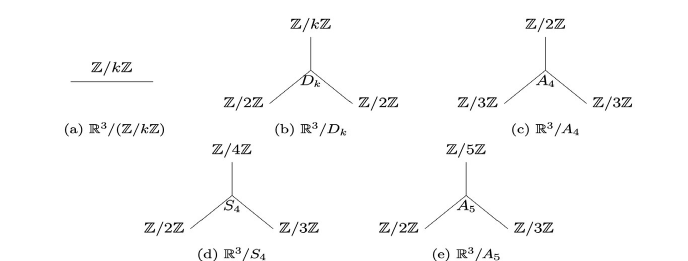}
\caption{Finite quotients of $\R^3$.}\label{fig:finitequotients}
\end{figure}

It will also be useful for us to describe the quotient of $\R^3$ under a finite subgroup $\Gamma$ of $SO(3)$. In Figure~\ref{fig:finitequotients}, a line segment represents a stratum of $\R^3/\Gamma$ with indicated cyclic isotropy, the origin has isotropy $\Gamma$ and the complement has trivial isotropy.

We also describe effective $O(2)$-representations on $\R^4$.

\begin{theorem}\label{thm:o2rep}
Let $\rho:O(2)\to O(4)$ be an effective representation of $O(2)$ and let $SO(2)\subset O(2)$ act as
\[
\rho(R(\theta))=\begin{bmatrix}R(p\theta)&0\\0&R(q\theta)\end{bmatrix},\qquad
R(\theta)=\begin{bmatrix}\cos\theta&-\sin\theta\\ \sin\theta&\cos\theta\end{bmatrix},
\]
where $p,q$ are coprime integers. Take $\tau\in O(2)\setminus SO(2)$ such that $\tau^2=\id$. Then, up to conjugation, changes of signs of $p,q$, and permutation of the two summands, $\rho$ lies in one of the following categories:
\begin{enumerate}[label=(\arabic*),leftmargin=2.2em]
\item $p,q\neq0$, and $\rho(\tau)=\diag(1,-1,1,-1)$;
\item $p=0$, $q=1$, and $\rho(\tau)=\diag(1,1,1,-1)$;
\item $p=0$, $q=1$, and $\rho(\tau)=\diag(1,-1,1,-1)$;
\item $p=0$, $q=1$, and $\rho(\tau)=\diag(-1,-1,1,-1)$.
\end{enumerate}
\end{theorem}

\begin{proof}
If both weights are nonzero, $\tau$ conjugates each rotation to its inverse. If the absolute values of the two weights are different, it preserves the two weight planes and acts as a reflection on each. Suppose that the weights have the same absolute value. After changing signs, the restriction of $SO(2)$ is scalar multiplication on $\C^2$, and $\tau$ is an orthogonal involution which anti-commutes with the complex structure. Its $+1$-eigenspace is a 2-dimensional totally real subspace $L$, and $\C^2=L\oplus iL$. Choosing an orthonormal basis of $L$ gives $\rho(\tau)=\diag(1,-1,1,-1)$. This gives (1).

Suppose now that $p=0$ and $q=1$. On the rotating plane, $\tau$ has to act as a reflection. On the zero-weight plane, it is an orthogonal involution, and hence its $+1$-eigenspace has dimension 2, 1, or 0. These are precisely (2), (3), and (4).
\end{proof}

\section{The Structure of Orbit Spaces and Orbit Types}

Throughout the remainder of the paper, $G$ will always be the Lie group $SU(2)$ or $SO(3)$, and $M$ is a simply connected closed 6-dimensional Riemannian manifold with positive sectional curvature which admits an effective isometric $G$-action.

We start by observing the following dichotomy for the topology of the orbit space $M/G$.

\begin{theorem}\label{thm:orbitspace}
The orbit space $M^*$ is homeomorphic to either $S^3$, a 3-ball $B^3$, or $B^4$. When $M^*=B^4$, $M$ is equivariantly diffeomorphic to $S^6$ with a fixed point homogeneous linear $SO(3)$-action.
\end{theorem}

\begin{proof}
The cohomogeneity is calculated via
\[
\dim(M^*)=\dim(M)-\dim(G)+\dim(H)=3+\dim(H),
\]
where $H$ is the principal isotropy group. $H$ is either 0- or 1-dimensional, since closed subgroups of $G=SU(2),SO(3)$ have dimensions 0, 1, 3, and $H$ cannot be 3-dimensional since otherwise the $G$-action would be trivial. Thus the cohomogeneity is either 3 or 4.

Suppose that the cohomogeneity is 4. Then the principal isotropy group $H$ is 1-dimensional, thus one of $S^1$, $O(2)$ or $\Pin(2)$. Since the isotropy representation of $G/H$ is irreducible, Theorem~\ref{thm:isotropy} implies that the isotropy representation of $G/H$ is equivalent to a subrepresentation of the isotropy representation of $K/H$, where $K\subset G$ is an isotropy group of a boundary face. The only possibility is $K=G$ and thus one boundary face has isotropy $K=G$, which means that the $G$-action is fixed point homogeneous. From Corollary~\ref{cor:lowfixed}, it follows that the only fixed point homogeneous action of $G$ on $M$ with cohomogeneity 4 is the linear $SO(3)$-action on $S^6$ described in Example~\ref{ex:s6low}(2).

When the cohomogeneity is 3, Theorem~\ref{thm:3manifold} implies that the orbit space $M^*$ is a simply connected 3-dimensional topological manifold possibly with boundary. If $\partial M^*$ is non-empty, let $s_0$ be the soul point. Since $s_0$ is an interior point of the topological 3-manifold $M^*$, its space of directions is homeomorphic to $S^2$. Theorem~\ref{thm:boundarysoul}(2) therefore gives $\partial M^*\cong S^2$. Capping $M^*$ by a 3-ball gives a closed simply connected 3-manifold, hence a 3-sphere by the Poincar\'e theorem. The 3-dimensional Schoenflies theorem now implies that $M^*$ is a 3-ball. If $\partial M^*$ is empty, then $M^*$ is a simply connected closed 3-manifold, thus a 3-sphere.
\end{proof}

We note that we have four kinds of orbits, corresponding to the 0-, 1-, or 3-dimensional closed subgroups of $G$:
\begin{enumerate}[label=(\arabic*),leftmargin=2.2em]
\item principal orbits $G/H$, with principal isotropy group $H$;
\item exceptional orbits $G/\Gamma$, with isotropy groups $\Gamma$ which are finite extensions of $H$;
\item singular orbits $G/K$, with 1-dimensional isotropy groups $K$, and hence $K=SO(2),O(2)$ when $G=SO(3)$, and $K=U(1),\Pin(2)$ when $G=SU(2)$;
\item fixed points, i.e., $G_p=G$.
\end{enumerate}

We start by showing that the principal isotropy is trivial apart from the two linear actions in Theorem~\ref{thm:main}(2).

\begin{theorem}\label{thm:principal}
If the principal isotropy subgroup $H$ is non-trivial, then $G=SO(3)$ and $M$ is equivariantly diffeomorphic to $S^6$ with one of the linear actions in Examples~\ref{ex:s6low}(2)--(3).
\end{theorem}

\begin{proof}
We break up the proof into several lemmas.

\begin{lemma}\label{lem:fixedprincipal}
If $M^G\neq\emptyset$ and the principal isotropy $H$ is non-trivial, then the $G$-action on $M$ has fixed point cohomogeneity at most one.
\end{lemma}

\begin{proof}
We separate the cases of $SU(2)$ and $SO(3)$ actions.

\emph{Case 1: $G=SU(2)$.} The action of $SU(2)$ on the normal space to $M^G$ is faithful and has no fixed vectors. By considering faithful real representations without trivial summands of $SU(2)$ in dimensions at most 6, we see that only the 4-dimensional irreducible representation of $SU(2)$ satisfies the requirements. This representation is the realification of the standard $SU(2)$-action on $\C^2$, and hence $M$ is fixed point homogeneous. But then $H$ is trivial, contradicting our assumption.

\emph{Case 2: $G=SO(3)$.} For the action of $SO(3)$ on the normal space to $M^G$ we have the following possibilities:
\begin{enumerate}[label=(\arabic*),leftmargin=2.2em]
\item $\R^3$ with the standard $SO(3)$-action. In this case the action is fixed point homogeneous;
\item $\R^5$ with the unique 5-dimensional irreducible representation of $SO(3)$. The $SO(3)$-action on the unit normal sphere $S^4$ has cohomogeneity one, which by definition implies that the $G$-action on $M$ has fixed point cohomogeneity one;
\item $\R^3\oplus\R^3$ with diagonal action of $SO(3)$. But in this case the principal isotropy is trivial, since $(x,y)\in\R^3\oplus\R^3$ has trivial isotropy when $x$ and $y$ are linearly independent.
\end{enumerate}
From Corollary~\ref{cor:lowfixed}, the only actions with fixed point cohomogeneity at most one and with non-trivial principal isotropy are the $SO(3)$-actions on $S^6$ described in Examples~\ref{ex:s6low}(2)--(3).
\end{proof}

In Lemmas~\ref{lem:nofixedprincipal} and \ref{lem:excludeprincipal} we will show that the case $M^G=\emptyset$ does not occur when the principal isotropy $H$ is non-trivial.

\begin{lemma}\label{lem:nofixedprincipal}
If $M^G=\emptyset$ and $H$ is non-trivial, then $G=SO(3)$, $H=\Z/2\Z$ or $\Z/2\Z\oplus\Z/2\Z$, and $\partial M^*$ consists of one face with isotropy $O(2)$. Furthermore, the interior has at most one singular orbit.
\end{lemma}

\begin{proof}
By Theorem~\ref{thm:orbitspace}, we only need to consider the case $\dim(M^*)=3$ and thus $H$ is finite. First notice that the Isotropy Lemma implies that $\partial M^*\neq\emptyset$, and thus $M^*=B^3$. If the isotropy representation of $H$ on the tangent space to $G/H$ has an irreducible subrepresentation of dimension greater than 1, then by the Isotropy Lemma the isotropy $K$ of a boundary face has dimension at least 2. Hence $K=G$, i.e., the $G$-action is fixed point homogeneous, which contradicts our assumption $M^G=\emptyset$.

In all other cases the irreducible components of the 3-dimensional isotropy representation of $G/H$ are 1-dimensional. For $G=SU(2)$, effectiveness implies that the center is not contained in $H$. Hence $H$ is cyclic of odd order, and the adjoint representation contains a nontrivial 2-dimensional rotation representation unless $H$ is trivial. Thus the $SU(2)$ case cannot occur. Among the non-trivial subgroups of $SO(3)$, only $\Z/2\Z$ or $\Z/2\Z\oplus\Z/2\Z$ has a 3-dimensional representation of this type. Thus $H=\Z/2\Z$ or $\Z/2\Z\oplus\Z/2\Z$, and $G=SO(3)$.

Since the isotropy representation of $SO(3)/H$ has a 1-dimensional subrepresentation on which $H$ acts as $-\id$, the Isotropy Lemma implies that the boundary face has isotropy $K=O(2)$, with $H$ meeting the non-identity component of $O(2)$. The only strata with higher-dimensional isotropy groups on $\partial M^*$ are fixed points, which cannot occur by our assumption. So the boundary has a single orbit type $O(2)$. From Theorem~\ref{thm:boundarysoul}(3), it follows that $M^*$ has at most one interior singular orbit, which would have to be the soul point.
\end{proof}

\begin{lemma}\label{lem:excludeprincipal}
If $M^G$ is empty and there are no interior singular orbits, then $M$ is not simply connected. Moreover, there are no actions with one interior singular orbit.
\end{lemma}

\begin{proof}
Recall that according to Lemma~\ref{lem:nofixedprincipal}, we have $G=SO(3)$. Let $\pi:M\to M^*$ be the natural projection, let $U=\pi^{-1}(\Int(M^*))$, and let $V$ be a tubular neighborhood of $\pi^{-1}(\partial M^*)$. Then $M=U\cup V$.

Suppose first that $M^*$ has no interior singular orbit. Then $U$ is an $SO(3)/H$-bundle over $\Int(M^*)=D^3$ and thus $U=SO(3)/H\times D^3$. The space $V$ deformation retracts onto an $SO(3)/O(2)$-bundle over $\partial M^*=S^2$, with structure group $N(O(2))/O(2)=\id$. Thus $V$ retracts onto the trivial $\RP^2$-bundle over $S^2$. Moreover, $U\cap V$ deformation retracts onto $SO(3)/H\times S^2$. Let $i:U\cap V\to U$ and $j:U\cap V\to V$ denote the respective inclusions. Since $i_*$ is an isomorphism on fundamental groups, van Kampen's theorem gives $\pi_1(M)\cong\pi_1(V)\cong\Z/2\Z$. This contradicts the assumption that $M$ is simply connected.

Suppose now that $M^*$ has an interior singular orbit. A priori the singular orbit could have isotropy $SO(2)$ or $O(2)$. Take the $SO(2)$-fixed point set. Over the boundary, $SO(2)$ fixes one point in every $SO(3)/O(2)$-orbit, so the corresponding fixed set is homeomorphic to $S^2$. If the interior singular orbit has isotropy $O(2)$, it contributes one further fixed point. Hence $\chi(M)=\chi(M^{SO(2)})=3$, contradicting Theorem~\ref{thm:eulerfixed}(2). Thus the interior singular orbit has isotropy $SO(2)$.

Suppose that the slice representation of $SO(2)$ has slope $(p,q)$, i.e., it acts on $\C\oplus\C$ as
\[
(z_1,z_2)\longmapsto(\xi^p z_1,\xi^q z_2),\qquad |\xi|=1.
\]
Then $H=\Z/\gcd(p,q)\Z$, and in particular $H$ is cyclic. Thus $H=\Z/2\Z$. The interior singular orbit contributes two $SO(2)$-fixed points, and therefore $\chi(M)=4$.

We first show that $M$ has the cohomology groups of $\CP^3$, and then use the Mayer--Vietoris sequence to obtain a contradiction. The fixed point set $M^H$ has dimension at least 3 since $H$ fixes a point in each principal orbit. Since $H$ is generated by an orientation-preserving involution, its fixed point components have even codimension. Hence the relevant component of $M^H$ is 4-dimensional. It is connected by Frankel's theorem \cite{Fra61}, and Theorem~\ref{thm:connectedness} implies that the inclusion $M^H\hookrightarrow M$ is 3-connected. In particular, $M^H$ is simply connected. Moreover, $N(H)/H=O(2)$ acts effectively on $M^H$, so $M^H$ is homeomorphic to $S^4$ or $\CP^2$ by \cite{HK89}. If $M^H$ is homeomorphic to $S^4$, then the connectedness lemma implies that $M$ is a homology 6-sphere, violating $\chi(M)=4$. Thus $M^H$ is homeomorphic to $\CP^2$. We then have $\pi_2(M)=\Z$ and $\pi_3(M)=0$. The Hurewicz theorem gives $H_2(M;\Z)=\Z$. Since $\chi(M)=4$, Poincar\'e duality and the universal coefficient theorem give
\[
H^0(M)=H^2(M)=H^4(M)=H^6(M)=\Z,\qquad H^{2i+1}(M)=0.
\]
Thus $M$ has the cohomology groups of $\CP^3$.

Now apply the Mayer--Vietoris sequence in integral cohomology. From the slice theorem,
\[
U=SO(3)\times_{SO(2)}D^4,
\]
so $U$ deformation retracts onto $S^2$. The space $V$ is homotopy equivalent to $S^2\times\RP^2$. The intersection $U\cap V$ is homotopy equivalent to $SO(3)\times_{SO(2)}S^3$, an $S^3$-bundle over $S^2$; in particular $H^2(U\cap V;\Z)\cong\Z$. The relevant part of the cohomology sequence is
\[
0\longrightarrow H^2(M;\Z)\longrightarrow H^2(U;\Z)\oplus H^2(V;\Z)
\longrightarrow H^2(U\cap V;\Z)\longrightarrow H^3(M;\Z)\longrightarrow0.
\]
Here
\[
H^2(M;\Z)\cong\Z,\qquad
H^2(U;\Z)\oplus H^2(V;\Z)\cong\Z\oplus\Z\oplus\Z/2\Z,
\qquad H^2(U\cap V;\Z)\cong\Z,
\]
and $H^3(M;\Z)=0$. The torsion summand maps trivially to the final copy of $\Z$, and therefore belongs to the kernel of the middle map. This kernel cannot be isomorphic to $\Z$, a contradiction.
\end{proof}

Combining Lemmas~\ref{lem:fixedprincipal}, \ref{lem:nofixedprincipal} and \ref{lem:excludeprincipal} finishes the proof of Theorem~\ref{thm:principal}.
\end{proof}

In our setting, exceptional orbits can have rich structure, making it difficult to recover the original manifold from the orbit space. We record the restrictions which hold without excluding the polyhedral subgroups of $SO(3)$.

\begin{proposition}\label{prop:exceptional}
Assume that the principal isotropy is trivial.
\begin{enumerate}[label=(\arabic*),leftmargin=2.2em]
\item If $G=SO(3)$, an exceptional isotropy group is one of $\Z/m\Z$, $D_m$, $A_4$, $S_4$, $A_5$. If $G=SU(2)$, an exceptional isotropy group is a finite subgroup of $SU(2)$; Proposition~\ref{prop:su2center} will show that it is cyclic of odd order.
\item In the orbit space $M/G$, every exceptional stratum has dimension at most one, and the closure of every connected component of the exceptional set contains a singular stratum.
\end{enumerate}
Moreover, if $C\subset\Gamma\subset G$ and the subgroups of $\Gamma$ which are $G$-conjugate to $C$ form one $\Gamma$-conjugacy class, then
\begin{equation}\label{eq:fixedhomogeneous}
(G/\Gamma)^C\cong N_G(C)/N_\Gamma(C).
\end{equation}
\end{proposition}

\begin{proof}
For $G=SO(3)$, the first assertion follows from the classification of finite subgroups listed in Section 2. For $G=SU(2)$, if $M^G\neq\emptyset$, Proposition~\ref{prop:fixedpoints} and a direct inspection of the two linear actions show that there are no exceptional orbits. If $M^G=\emptyset$, Proposition~\ref{prop:su2center} below shows that an exceptional isotropy group does not contain the center and hence is cyclic of odd order. Notice that Proposition~\ref{prop:su2center} does not use the present proposition. Formula~\eqref{eq:fixedhomogeneous} follows directly from
\[
(G/\Gamma)^C=\{g\Gamma\mid g^{-1}Cg\subset\Gamma\}.
\]
For example, for the polyhedral groups in $SO(3)$ one has
\[
N_{A_4}(C_2)=\Z/2\Z\oplus\Z/2\Z,\qquad
N_{S_4}(C_3)=D_3,\qquad
N_{A_5}(C_5)=D_5.
\]
Since $N_{SO(3)}(C_k)=O(2)$, the corresponding fixed sets in $SO(3)/A_4$, $SO(3)/S_4$ and $SO(3)/A_5$ are single circles.

By Theorem~\ref{thm:nospecial}, exceptional strata cannot have codimension one in the 3-dimensional orbit space. Hence they have dimension at most one. At an exceptional orbit $G/\Gamma$, the slice is 3-dimensional. The action of $\Gamma$ on the tangent space of the orbit is orientation-preserving, and the action on $M$ is orientation-preserving since $G$ is connected. Thus the slice representation has image in $SO(3)$. Every non-trivial finite subgroup of $SO(3)$ has a rotation axis, so an exceptional point is contained in the closure of a 1-dimensional exceptional stratum. Consequently the exceptional set is a finite graph with no isolated vertices.

Suppose that one connected component $E$ of the exceptional set has closure disjoint from all singular strata. In the $SU(2)$ case, the preceding paragraph and Proposition~\ref{prop:su2center} allow us to choose a non-central cyclic isotropy subgroup. Choose such a non-trivial cyclic subgroup $C$ which is the isotropy group along an open edge of $E$, and let $F$ be a component of $M^C$ which projects to that edge. Along the edge, $(G/\Gamma)^C$ is a union of circles, and hence $F$ is 2-dimensional. The component $F$ is compact because the closure of $E$ contains no singular orbit. Moreover, the connected group $N_G(C)^0\cong S^1$ preserves $F$ and acts locally freely on it; a fixed point of this circle action would have positive-dimensional isotropy and hence would lie over a singular orbit. The infinitesimal action therefore gives a nowhere-zero vector field on the compact totally geodesic positively curved surface $F$. This gives $\chi(F)=0$, contradicting the Gauss--Bonnet theorem. Thus the closure of every exceptional component contains a singular stratum.
\end{proof}

The normalizer calculation in Proposition~\ref{prop:exceptional} shows that the fixed surface argument does not exclude $A_4$, $S_4$ or $A_5$. In fact, Example~\ref{ex:s6remaining}(3) below gives an exceptional $A_4$-orbit on the round $6$-sphere.

\section{Actions, Orbit Spaces and the Topology of $G$-manifolds}

In this section we study different types of $G$-actions on positively curved 6-manifolds. We start with the case of non-empty fixed point set.

\begin{proposition}\label{prop:fixedpoints}
If $M^G\neq\emptyset$, then one of the following holds:
\begin{enumerate}[label=(\arabic*),leftmargin=2.2em]
\item $M$ is equivariantly diffeomorphic to $S^6$ or $\CP^3$ with a linear action;
\item $G=SO(3)$ and $M^G$ is finite. In this case $M^*=B^3$ and $M^G$ lies on $\partial M^*$.
\end{enumerate}
\end{proposition}

\begin{proof}
We first note that a substantial part of the proof has already appeared in the proof of Lemma~\ref{lem:fixedprincipal}. In fact, we proved that if $G=SU(2)$, then $M$ is fixed point homogeneous. If $G=SO(3)$, either $M$ has fixed point cohomogeneity at most one, which is classified by Theorem~\ref{thm:fixedpointcohom}, or the action on the tangent space at a fixed point is given by
\[
A\cdot(x,y)=(Ax,Ay),\qquad A\in SO(3),\quad (x,y)\in\R^3\oplus\R^3.
\]
Thus the fixed points are isolated and the orbit types of this action are:
\begin{itemize}[leftmargin=2em]
\item principal orbits with trivial isotropy, represented by two linearly independent vectors in $\R^3$;
\item singular orbits with $SO(2)$-isotropy, represented by two linearly dependent vectors in $\R^3$ which are not both zero;
\item the fixed point $(0,0)$.
\end{itemize}
The union of singular orbits near a fixed point has dimension 4 in $M$, which descends to a 2-dimensional stratum of $M^*$. Since $\dim(M^*)=3$, this stratum is a boundary face in $\partial M^*$ and hence $M^*=B^3$ by Theorem~\ref{thm:orbitspace}. Moreover, all fixed points are contained in the boundary strata.
\end{proof}

Next, we study the case where $M^G$ is finite or empty.

\subsection{$SO(3)$ Actions with $M^G$ Finite or Empty}

By Theorem~\ref{thm:orbitspace}, we divide this subsection into two parts, corresponding to $M^*=B^3$ and $M^*=S^3$, and we start with the case where $M^*=B^3$.

\begin{theorem}\label{thm:so3ball}
Assume $G=SO(3)$, with $M^*=B^3$ and $M^G$ finite. Then:
\begin{enumerate}[label=(\arabic*),leftmargin=2.2em]
\item The principal isotropy is trivial and the boundary face of $M^*$ consists of singular orbits with $SO(2)$-isotropy.
\item If $\partial M^*$ contains more than one orbit type, then $\partial M^*$ contains exactly two singular points which are either two fixed points or one fixed point and one $O(2)$-orbit. Furthermore, these two singular orbits do not lie in a face.
\item If there is an $O(2)$-orbit on $\partial M^*$, then there exist an interior singular orbit and a $\Z/2\Z$-exceptional stratum connecting the interior singular orbit and the $O(2)$-orbit.
\item There is at most one interior singular orbit, whose isotropy group has to be $SO(2)$.
\item If $\partial M^*$ contains only one orbit type, then $\chi(M)=4$ or $6$. If $\partial M^*$ contains more than one orbit type, then $\chi(M)=2$ or $4$; in the case of an $O(2)$-orbit on the boundary, one has $\chi(M)=4$.
\end{enumerate}
\end{theorem}

\begin{proof}
\emph{Proof of (1).} If the principal isotropy is non-trivial, then by Theorem~\ref{thm:principal} and Corollary~\ref{cor:lowfixed}, the action on $M$ is given by Example~\ref{ex:s6low}(3). But then $M^G$ is a circle, contradicting our assumption that $M^G$ is finite. By Theorem~\ref{thm:nospecial}, the boundary does not contain exceptional orbits.

Suppose that the boundary face orbits have $O(2)$-isotropy. The slice action on the 4-dimensional normal space is either ineffective, in which case the principal isotropy is non-trivial, or effective. Since the $O(2)$-stratum is 2-dimensional, the effective slice action is case (2) of Theorem~\ref{thm:o2rep}. A generic point in the rotating plane is fixed by a reflection, and hence the principal isotropy contains $\Z/2\Z$, again a contradiction. Therefore the boundary face has $SO(2)$-isotropy.

\emph{Proof of (4).} Theorem~\ref{thm:boundarysoul}(3) implies that $\Int(M^*)$ has at most one isolated stratum, which is the soul point. By Proposition~\ref{prop:exceptional}(2), the soul point cannot be an exceptional orbit. Thus it has to be a singular orbit $G/K$. A priori $K$ could be $SO(2)$ or $O(2)$. Suppose $K=O(2)$. The slice representation of $O(2)$ is 4-dimensional and orientation-reversing, since the isotropy representation of $O(2)$ on $SO(3)/O(2)$ is orientation-reversing. We apply Theorem~\ref{thm:o2rep}.

In case (1), the involution acts by $\diag(1,-1,1,-1)$ and preserves the orientation of the slice. Case (3) is also orientation-preserving. In case (2), the $O(2)$-stratum is 2-dimensional and hence lies in $\partial M^*$. In case (4), the $SO(2)$-stratum is 2-dimensional and lies in $\partial M^*$, and the $O(2)$-orbit belongs to its closure. Thus none of the four cases gives an interior $O(2)$-orbit. In conclusion $K=SO(2)$.

\emph{Proof of (3).} Suppose there is an $O(2)$-orbit on $\partial M^*$. Its slice action is effective since otherwise the principal isotropy is non-trivial. The boundary face has isotropy $SO(2)$. Case (2) of Theorem~\ref{thm:o2rep} gives a 2-dimensional $O(2)$-stratum, while case (1) has no nearby $SO(2)$-strata. Case (3) has nearby $SO(2)$-strata, but it also has a codimension one exceptional $\Z/2\Z$-stratum: in coordinates $(a,b,z)\in\R\oplus\R\oplus\C$, this stratum is given by $b=0$ and $z\neq0$. This is a special exceptional stratum and is excluded by Theorem~\ref{thm:nospecial}. Hence the slice action is case (4). In particular, the $O(2)$-orbit is isolated and a $\Z/2\Z$-stratum emanates from it.

By Theorem~\ref{thm:boundarysoul}(3), this stratum either ends at the soul point or ends at another boundary stratum point. It cannot end at a fixed point, since the slice representation $\R^3\oplus\R^3$ at an isolated fixed point has no $\Z/2\Z$-strata.

We claim that there cannot be two $O(2)$-orbits on $\partial M^*$. Assume otherwise. First suppose that a $\Z/2\Z$-stratum connects the two $O(2)$-orbits. Let $C=\Z/2\Z$ be its isotropy group. We have
\[
(SO(3)/SO(2))^C=2\ \mathrm{pt},\qquad
(SO(3)/O(2))^C=S^1\cup\mathrm{pt},\qquad
(SO(3)/C)^C=S^1\cup S^1.
\]
Let $F$ be a component of $M^C$ containing one of the circle components over the exceptional stratum. The component $F$ is a compact 2-dimensional totally geodesic submanifold with positive curvature. The identity component $N(C)^0$ acts locally freely on $F$, including at the two end orbits where $F$ meets the circle component of $(SO(3)/O(2))^C$. Hence $F$ admits a nowhere-zero vector field, contradicting the Gauss--Bonnet theorem.

The remaining possibility is that two $\Z/2\Z$-strata connect the two $O(2)$-orbits to the soul point. By (4), the soul point has isotropy $SO(2)$. If the slice weights at the soul are $(p,q)$, the two distinct $\Z/2\Z$-strata force both $p$ and $q$ to be even. This forces the principal isotropy to contain $\Z/2\Z$, a contradiction. Therefore there is at most one $O(2)$-orbit on the boundary. Its exceptional stratum cannot end at a fixed point or at another $O(2)$-orbit, and hence it ends at the soul. This proves (3).

\emph{Proof of (2).} In each boundary face orbit, $SO(2)$ fixes exactly two points. In an $O(2)$-orbit or a $G$-fixed point, $SO(2)$ fixes one point. Thus the $SO(2)$-fixed point component over $\partial M^*$ is a branched double cover of $\partial M^*=S^2$, with branch points corresponding to $O(2)$-orbits or $G$-fixed points. If the boundary contains more than one orbit type, this cover has a branch point and is connected. It is an orientable positively curved surface and hence a 2-sphere. The Riemann--Hurwitz formula shows that a branched double cover between two 2-spheres has exactly two branch points. We have shown in (3) that there cannot be two $O(2)$-orbits. Thus the two branch points correspond to two fixed points or one fixed point and one $O(2)$-orbit.

\emph{Proof of (5).} If $\partial M^*$ has more than one orbit type, the component of $M^{SO(2)}$ over the boundary is a 2-sphere by the proof of (2). The only possible further contribution is the interior $SO(2)$-orbit, which contributes two points. Hence $\chi(M)=2$ or $4$. If an $O(2)$-orbit occurs on the boundary, (3) gives an interior $SO(2)$-orbit, and therefore $\chi(M)=4$.

If $\partial M^*$ has only one orbit type, the subset of $M^{SO(2)}$ over the boundary is an unbranched double cover of $S^2$, and hence is the union of two 2-spheres. The Boundary Soul Lemma shows that the only possible additional singular stratum is the soul point, which has $SO(2)$-isotropy by (4). Thus
\[
M^{SO(2)}=S^2\cup S^2\quad\text{or}\quad
M^{SO(2)}=S^2\cup S^2\cup\{2\ \mathrm{points}\},
\]
and $\chi(M)=4$ or $6$.
\end{proof}

\begin{figure}[ht]
\centering
\includegraphics[width=.70\textwidth]{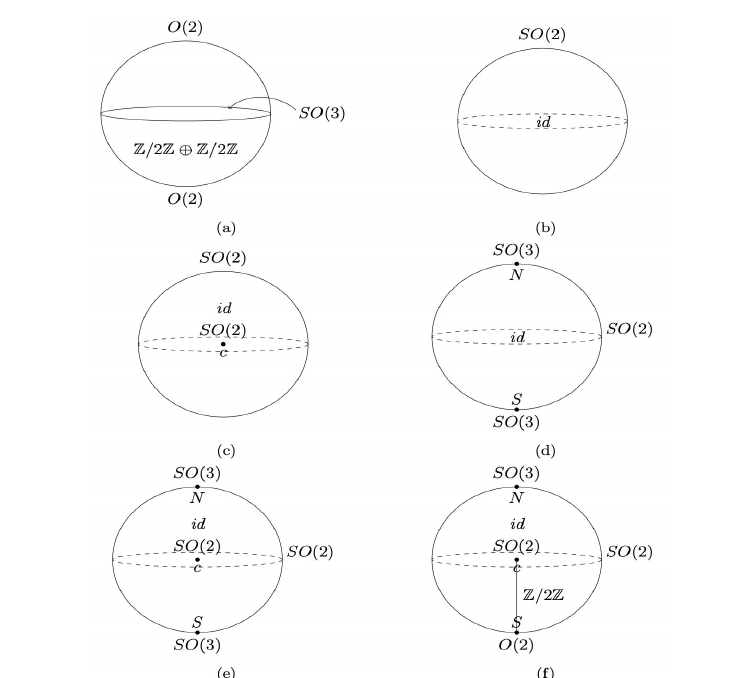}
\caption{$G=SO(3)$, $M^*=B^3$.}\label{fig:so3ball}
\end{figure}

\begin{corollary}\label{cor:so3ballpictures}
If $G=SO(3)$ and $M^*=B^3$, then the singular stratification, together with the exceptional stratum forced by an $O(2)$-orbit, is as in Figures~\ref{fig:so3ball}(a)--(f). Any further exceptional component has closure meeting one of the displayed singular strata. In the pictures, the groups represent the isotropy groups of the corresponding strata.
\end{corollary}

\begin{proof}
If $M^G$ is not finite, Proposition~\ref{prop:fixedpoints} implies that $M^*$ is as in Figure~\ref{fig:so3ball}(a). So we may assume that $M^G$ is finite. Then Theorem~\ref{thm:so3ball} implies that the boundary face has isotropy $SO(2)$ and that there is at most one interior singular orbit.

If $\partial M^*$ has one orbit type and $\Int(M^*)$ contains no singular orbit, the singular stratification is as in Figure~\ref{fig:so3ball}(b). If $\Int(M^*)$ contains one singular orbit, it is as in Figure~\ref{fig:so3ball}(c). If $\partial M^*$ has multiple orbit types, Theorem~\ref{thm:so3ball} implies that it contains two fixed points, or one fixed point and one $O(2)$-orbit. If it contains two fixed points, the singular stratification is as in Figure~\ref{fig:so3ball}(d) or (e). Otherwise, the forced strata are as in Figure~\ref{fig:so3ball}(f). The last assertion follows from Proposition~\ref{prop:exceptional}(2).
\end{proof}

We point out that in Figures~\ref{fig:so3ball}(a) and (d), $M$ is classified, see Proposition~\ref{prop:fixedpoints} and the following theorem respectively. We now prove Theorem~\ref{thm:so3main}.

\begin{theorem}\label{thm:so3classification}
Assume that $G=SO(3)$, $M^*=B^3$, $\partial M^*$ contains more than one orbit type and that there are no exceptional orbits or interior singular orbits. Then $M^6$ is equivariantly homeomorphic to $S^6$ with the linear $SO(3)$-action in Example~\ref{ex:s6remaining}(2).
\end{theorem}

\begin{proof}
Theorem~\ref{thm:so3ball}(2) implies that $\partial M^*$ has two singular points. Since $M$ has no exceptional orbits, these two orbits cannot be $O(2)$-orbits, because case (4) of Theorem~\ref{thm:o2rep} gives a nearby $\Z/2\Z$-stratum. Thus the two singular points are two $G$-fixed points.

The orbit types are now $(H)=(\id)$, $(K)=(SO(2))$ and $(G)=(SO(3))$, with two fixed points. We use Corollary V.6.2 and Proposition V.10.1 of \cite{Bre72}. After deleting the two fixed points, the classification reduces to $G$-spaces with orbit space a 2-disk, principal isotropy $H=\id$, and singular isotropy $K=SO(2)$ on the boundary. Their equivariant homeomorphism classes are parametrized by
\[
[\partial D^2,N(H)/(N(H)\cap N(K))]/\pi_0(N(H)/H)
=\pi_1(\RP^2)=\Z/2\Z.
\]
The two corresponding 5-dimensional $G$-spaces are the 5-sphere with the diagonal action on $\R^3\oplus\R^3$, and $S^2\times S^3$ with the other class. The original space $M$ is the suspension of one of these spaces. The suspension is a manifold only in the first case, since the link at a suspension point has to be a 5-sphere. Hence $M$ is equivariantly homeomorphic to the suspension of the diagonal action on $S^5$, which is the action in Example~\ref{ex:s6remaining}(2).
\end{proof}

\begin{remark}\label{rem:so3open}
The other cases of $G=SO(3)$ and $M^*=B^3$ listed in Corollary~\ref{cor:so3ballpictures} are also interesting but not yet fully understood. For Figure~\ref{fig:so3ball}(b), $\partial M^*$ has isotropy $SO(2)$ and $\Int(M^*)$ has trivial isotropy. Corollary V.6.2 in \cite{Bre72} implies that such $G$-spaces are parametrized by
\[
[B,N(H)/(N(H)\cap N(K))]/\pi_0(N(H)/H)=\pi_2(\RP^2)=\Z.
\]
Two examples are $\CP^3$ and $S^2\times S^4$. The $SO(3)$-action on $\CP^3$ is the one described in Example~\ref{ex:cp3}(2), while the action on $S^2\times S^4$ is diagonal with the standard $SO(3)$-action on the $S^2$-factor and the double suspension on the $S^4$-factor. We suspect the $G$-spaces are oriented $S^2$-bundles over $S^4$, which are classified by the first Pontryagin class. We do not know whether $S^2\times S^4$ with this $SO(3)$-action admits an invariant metric with positive sectional curvature. According to the Hopf conjecture, it should not have any metric with positive curvature.

For Figure~\ref{fig:so3ball}(f), Example~\ref{ex:cp3}(3) is an action with the corresponding orbit stratification. For Figures~\ref{fig:so3ball}(c) and (e), we do not have examples, and we expect that they do not exist.
\end{remark}

We now turn to the case of $G=SO(3)$ and $M^*=S^3$.

\begin{proposition}\label{prop:so3sphere}
Suppose that $G=SO(3)$ and $M^*=S^3$. Let $a$ be the number of singular orbits with $SO(2)$-isotropy and let $b$ be the number of singular orbits with $O(2)$-isotropy. Then
\[
a+b\leq3,\qquad \chi(M)=2a+b.
\]
In particular, $b$ is even. Thus $b=0$ or $2$, and in the latter case $a\leq1$.
\end{proposition}

\begin{proof}
There exist singular orbits since $\chi(M^{SO(2)})=\chi(M)>0$. Suppose there were four singular orbits. At a singular orbit $G/K_i$, the space of directions is $S^3(1)/K_i$, where $K_i=SO(2)$ or $O(2)$ acts linearly on the unit normal sphere. Lemma 4 in \cite{HK89} implies
\[
\xt_3(S^3/K_i)\leq \xt_3(S^2(1/2))=\frac\pi3.
\]
The same estimate for $O(2)$ follows also by passing from the $SO(2)$-quotient to the further quotient. Four singular orbits would contradict Lemma~\ref{lem:extent}. Hence $a+b\leq3$.

There are no fixed points because, for the representations which occur at an $SO(3)$-fixed point in dimension six, the image of the fixed point lies in the boundary of the orbit space. An $SO(2)$-orbit contributes two points to $M^{SO(2)}$, while an $O(2)$-orbit contributes one point. No principal or exceptional orbit has an $SO(2)$-fixed point. Therefore
\[
\chi(M)=\chi(M^{SO(2)})=2a+b.
\]
\end{proof}

\begin{remark}
We note that Fang already observed the bound on the number of singular orbits in Lemma 2.3 of \cite{Fang02}. We decide to include the proof here for the completeness of this paper.
\end{remark}

\begin{figure}[ht]
\centering
\includegraphics[width=.43\textwidth]{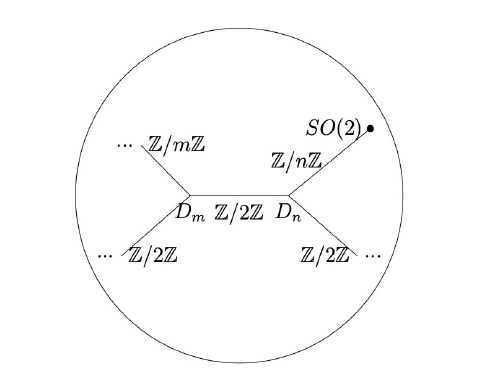}
\caption{$G=SO(3)$, $M^*=S^3$.}\label{fig:so3sphere}
\end{figure}

\begin{remark}
When $G=SO(3)$ and $M^*=S^3$, the possible stratification of $M^*$ is depicted in Figure~\ref{fig:so3sphere}. In the picture, $m,n$ are positive integers. Examples~\ref{ex:cp3}(4) and \ref{ex:flag}(2) are actions of this type.
\end{remark}

\subsection{$SU(2)$ Actions with $M^G=\emptyset$}

Let $Z=\{\pm1\}$ denote the center of $SU(2)$.

\begin{proposition}\label{prop:su2center}
Suppose that $G=SU(2)$ and $M^G=\emptyset$. Then:
\begin{enumerate}[label=(\arabic*),leftmargin=2.2em]
\item $M^*=S^3$.
\item The fixed point set $M^Z$ is precisely the union of the singular orbits. Every singular orbit has isotropy $U(1)$ and is a 2-sphere.
\item Every exceptional isotropy group is cyclic of odd order.
\item If $r$ is the number of singular orbits, then
\[
1\leq r\leq3,\qquad \chi(M)=2r.
\]
\end{enumerate}
\end{proposition}

\begin{proof}
By Theorem~\ref{thm:principal}, the principal isotropy is trivial. Since the action of $Z$ preserves orientation, every component of $M^Z$ has even codimension. A 0-dimensional component would be preserved pointwise by the connected group $SU(2)$ and hence would be contained in $M^G$. A 6-dimensional component would make the action ineffective. Thus the only possible dimensions are 2 and 4.

We first rule out a 4-dimensional component $F$. Since $F$ has codimension two, Theorem~\ref{thm:connectedness} implies that $F$ is simply connected. The induced action of $SO(3)=SU(2)/Z$ on $F$ is effective: its kernel is a normal subgroup of $SO(3)$, and the whole group cannot fix $F$ because $M^G=\emptyset$. Its cohomogeneity is one or two. By \cite{HK89}, $F$ is homeomorphic to $S^4$ or $\CP^2$.

We will use two elementary facts about $\Pin(2)$. First, there is no 2-dimensional real representation
\[
\rho:\Pin(2)\longrightarrow O(2)
\]
for which $\rho(-1)=-I$. Indeed, on $U(1)\subset\Pin(2)$ such a representation has the form $\rho(e^{i\theta})=R(m\theta)$ with $m$ odd. If $j\in\Pin(2)\setminus U(1)$, then $je^{i\theta}j^{-1}=e^{-i\theta}$, so $\rho(j)$ is a reflection and $\rho(j)^2=I$. On the other hand $j^2=-1$, which would give $\rho(j)^2=-I$, a contradiction.

Second, if $W$ is a 4-dimensional real representation of $\Pin(2)$ on which $-1$ acts as $-I$, then every $j\in\Pin(2)\setminus U(1)$ preserves the orientation of $W$. To see this, use the $U(1)$-action to identify $W$ with $\C^2$. The element $j$ is anti-complex-linear, and an anti-complex-linear automorphism of $\C^2$ has positive real determinant.

Suppose first that the induced $SO(3)$-action on $F$ has cohomogeneity one. Its two singular isotropy groups are $SO(2)$ or $O(2)$. By the double disk-bundle decomposition, if they were both $SO(2)$, then
\[
\chi(F)=\chi(SO(3)/SO(2))+\chi(SO(3)/SO(2))=4,
\]
contradicting $F\cong S^4$ or $\CP^2$. Hence one singular isotropy group is $O(2)$. It lifts to $\Pin(2)$ in $SU(2)$, while $Z$ acts as $-I$ on the normal 2-plane to $F$. This contradicts the observation above.

Suppose next that the induced action has cohomogeneity two. Its principal isotropy has dimension one. An $O(2)$-orbit gives the same contradiction, and an $SO(3)$-fixed point would be a $G$-fixed point in $M$. There are therefore no exceptional or singular orbit types, and every orbit has isotropy $SO(2)$. The orbit map is an $S^2$-bundle over a closed surface, with structure group $N(SO(2))/SO(2)=\Z/2\Z$. Since $F$ is simply connected, the base is $S^2$; since the structure group is discrete, the bundle is trivial. Hence $\chi(F)=4$, again a contradiction. We conclude that every component of $M^Z$ is 2-dimensional.

Let $F$ now be a 2-dimensional component of $M^Z$. For $x\in F$, the orbit $G\cdot x$ is contained in $F$. It cannot be 3-dimensional, so $G_x$ is 1-dimensional. Since $M^G=\emptyset$, it follows that $G_x=U(1)$ or $\Pin(2)$ and $F=G\cdot x$.

The second possibility cannot occur. If $G_x=\Pin(2)$, choose $j\in\Pin(2)\setminus U(1)$. The action of $j$ on the tangent plane of $SU(2)/\Pin(2)=\RP^2$ reverses orientation. On the 4-dimensional normal slice, the center acts as $-I$, and the second observation above shows that $j$ preserves orientation. Hence $j$ reverses the orientation of $T_xM$, contradicting the fact that every element of the connected group $SU(2)$ acts orientation-preservingly on $M$. Thus every singular isotropy group is $U(1)$ and every singular orbit is $S^2$.

It follows that $M^Z$ is precisely the union of all singular orbits. Moreover, if an exceptional isotropy group contained $Z$, the corresponding 3-dimensional orbit would be contained in the 2-dimensional set $M^Z$. This is impossible. A finite subgroup of $SU(2)$ not containing $Z$ is cyclic of odd order, which proves (2) and (3).

We now show that $M^*=S^3$. If $M^*=B^3$, a boundary face has positive-dimensional isotropy. Its preimage is 4-dimensional and is contained in $M^Z$, contradicting what we have just proved. Thus $M^*$ has no boundary and Theorem~\ref{thm:orbitspace} gives $M^*=S^3$.

There is at least one singular orbit because $\chi(M^{U(1)})=\chi(M)>0$. The Extent Lemma gives at most three singular orbits exactly as in Proposition~\ref{prop:so3sphere}. Each $U(1)$-orbit contains two $U(1)$-fixed points, and no other orbit contributes. Hence, if $r$ is the number of singular orbits,
\[
\chi(M)=\chi(M^{U(1)})=2r.
\]
\end{proof}

\begin{figure}[ht]
\centering
\includegraphics[width=.68\textwidth]{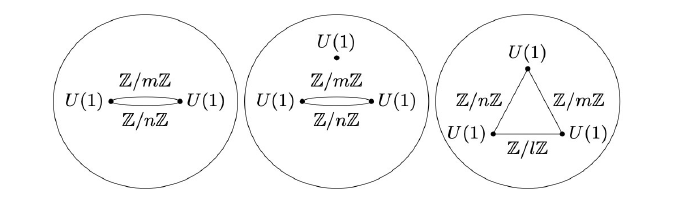}
\caption{$G=SU(2)$, $M^*=S^3$, with $U(1)$ singular isotropy.}\label{fig:su2sphere}
\end{figure}

\begin{remark}
When $G=SU(2)$ and $M^*=S^3$, the singular isotropy groups are all $U(1)$, and the possible stratification is depicted in Figure~\ref{fig:su2sphere}. Here $m,n,l$ are pairwise coprime odd integers. Example~\ref{ex:biquotient} is an action of this type.
\end{remark}

We now prove Theorem~\ref{thm:su2main}.

\begin{theorem}\label{thm:su2classification}
Let $G=SU(2)$.
\begin{enumerate}[label=(\arabic*),leftmargin=2.2em]
\item If the fixed point set $M^G$ is non-empty, then $M$ is equivariantly diffeomorphic to a linear action on $S^6$ or $\CP^3$.
\item If $M^G$ is empty and the action has no exceptional orbits, then the action is equivariantly diffeomorphic to one of the actions on $S^6$, $S^2\times S^4$ or $SU(3)/T^2$ described in Remark~\ref{rem:su2models}.
\end{enumerate}
\end{theorem}

\begin{proof}
Part (1) has already been proved in Proposition~\ref{prop:fixedpoints}. For (2), Proposition~\ref{prop:su2center} gives $M^*=S^3$ and between one and three singular orbits, all with $U(1)$-isotropy. At a singular orbit, write the slice representation as
\[
z\cdot(u,v)=(z^p u,z^q v),\qquad (u,v)\in\C^2.
\]
The singular orbit is isolated, so $p,q\neq0$. If $|p|>1$ or $|q|>1$, a coordinate axis gives an exceptional orbit. Hence $|p|=|q|=1$, and as a real representation the slice is the Hopf representation.

Let $r$ be the number of singular orbits. By Proposition~\ref{prop:su2center}, $1\leq r\leq3$ and $\chi(M)=2r$. A tubular neighborhood of each singular orbit is
\[
V=SU(2)\times_{U(1)}D^4.
\]
The associated oriented 4-plane bundle over $SU(2)/U(1)=S^2$ is trivial: the two weight-one planes give the zero class in $\pi_1(SO(4))=\Z/2\Z$. Thus $V$ is diffeomorphic to $S^2\times D^4$. Its boundary is a free $SU(2)$-manifold over $S^2$, and the corresponding principal $SU(2)$-bundle is trivial because $\pi_1(SU(2))=0$. Hence
\[
\partial V\cong SU(2)\times S^2
\]
equivariantly.

Remove disjoint tubular neighborhoods of the singular orbits. The remaining principal part is a principal $SU(2)$-bundle over $X_r=S^3\setminus\bigcup_{i=1}^r\Int(B_i^3)$. Since $X_r$ has the homotopy type of a wedge of $(r-1)$ copies of $S^2$, and $\pi_1(BSU(2))=\pi_2(BSU(2))=0$, this bundle is trivial. Therefore the manifold is obtained by equivariantly gluing $r$ copies of $V$ to $SU(2)\times X_r$.

An equivariant diffeomorphism of a boundary component has the form
\[
f(g,x)=(g\,a(x),\phi(x)),\qquad a:S^2\to SU(2),\quad \phi\in\Diff(S^2).
\]
The map $a$ is null-homotopic, and an orientation-preserving $\phi$ is isotopic to the identity by \cite{Sma59}. These isotopies extend equivariantly over a collar. An orientation-reversing representative also extends over $V$. Indeed, choose $j\in\Pin(2)\setminus U(1)$ and let $c:\C^2\to\C^2$ be complex conjugation. Then
\[
[g,v]\longmapsto[gj,c(v)]
\]
is a well-defined equivariant diffeomorphism of $SU(2)\times_{U(1)}D^4$ and induces an orientation-reversing diffeomorphism on $S^3/U(1)=S^2$. Thus all attaching maps are equivalent.

If $r=1$, the resulting decomposition is
\[
S^3\times D^3\ \cup_{S^3\times S^2}\ D^4\times S^2=\partial(D^4\times D^3)=S^6.
\]
If $r=2$, it is the double of $S^2\times D^4$, and hence $S^2\times S^4$. If $r=3$, there is again a unique equivariant diffeomorphism type. Example~\ref{ex:flag}(1) realizes it on $SU(3)/T^2$. This proves the theorem.
\end{proof}

\begin{remark}\label{rem:su2models}
The strategy of our proof is motivated by the proof of Theorem B in \cite{Fang02}. The $SU(2)$-actions on $S^6$, $S^2\times S^4$ and $SU(3)/T^2$ in the three cases are all realizable. On $S^6$ it is the triple suspension of the Hopf action on $S^3$. On $S^2\times S^4$ it is the diagonal action where $SU(2)$ acts as $SO(3)$ on the $S^2$-factor and acts on $S^4$ as the suspension of the Hopf action. On $SU(3)/T^2$ it acts via left multiplication. We do not know whether $S^2\times S^4$ admits a metric with positive sectional curvature invariant under this action. The actions on $S^6$ and $SU(3)/T^2$ preserve positive curvature.
\end{remark}

We now summarize the claim about the Euler characteristic.

\begin{theorem}\label{thm:euler}
Let $M=M^6$ be a 6-dimensional closed simply connected Riemannian manifold of positive sectional curvature such that $SU(2)$ or $SO(3)$ acts isometrically and effectively on $M$. Then the Euler characteristic $\chi(M)=2,4$, or $6$.
\end{theorem}

\begin{proof}
First recall that $M^*=B^3,S^3$ or $B^4$ by Theorem~\ref{thm:orbitspace}. The case $M^*=B^4$ is the fixed point homogeneous $SO(3)$-action on $S^6$. If $M^*=B^3$ and $M^G$ is not finite, Proposition~\ref{prop:fixedpoints} gives a linear action on $S^6$ or $\CP^3$. If $M^*=B^3$ and $M^G$ is finite, Theorem~\ref{thm:so3ball}(5) gives $\chi(M)\leq6$. Finally, when $M^*=S^3$, Propositions~\ref{prop:so3sphere} and \ref{prop:su2center} give $\chi(M)\leq6$. By Theorem~\ref{thm:eulerfixed}(2), $\chi(M)$ is positive and even. Hence $\chi(M)=2,4$, or $6$.
\end{proof}

Theorem~\ref{thm:principal}, Proposition~\ref{prop:exceptional}, Proposition~\ref{prop:su2center}, Proposition~\ref{prop:so3sphere}, and Theorem~\ref{thm:euler} prove Theorem~\ref{thm:main}. Theorems~\ref{thm:su2classification} and \ref{thm:so3classification} prove Theorems~\ref{thm:su2main} and \ref{thm:so3main}. Collecting all the results, we see that the remaining singular configurations are described in Figures~\ref{fig:so3ball}--\ref{fig:su2sphere}.

\section{Explicit Examples of $G$-actions}\label{sec:examples}

In this section we list the relevant isometric $SU(2)$- and $SO(3)$-actions on the known positively curved 6-manifolds, namely $S^6$, $\CP^3$, $SU(3)/T^2$ and $SU(3)//T^2$, and depict their orbit stratifications. For $S^6$ and $\CP^3$ we list all linear actions. For the known positively curved metrics on $SU(3)/T^2$ and $SU(3)//T^2$, the full isometry group was determined in \cite{GSZ06}, and one easily sees that the only isometric actions are the ones described below.

\begin{figure}[ht]
\centering
\includegraphics[width=.72\textwidth]{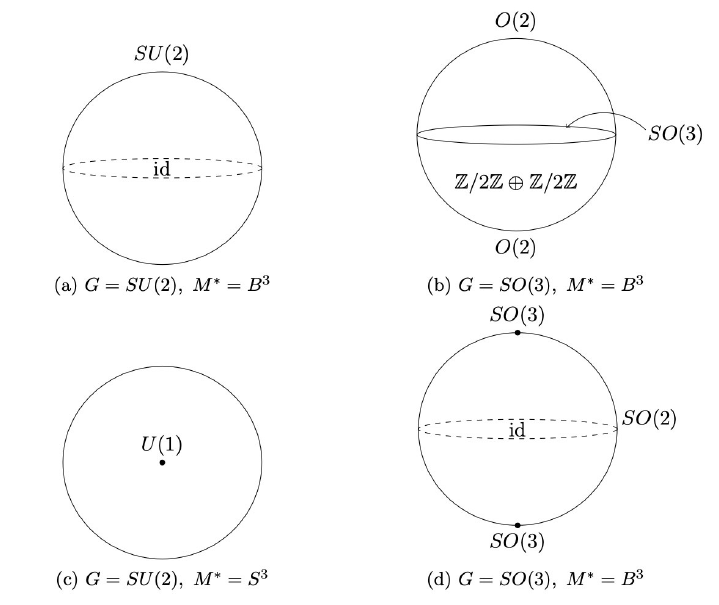}
\caption{$M=S^6$ for the four actions in Examples~\ref{ex:s6low}(1), \ref{ex:s6low}(3), \ref{ex:s6remaining}(1), and \ref{ex:s6remaining}(2).}\label{fig:s6actions}
\end{figure}

\begin{example}\label{ex:s6low}
Actions on $S^6$ with fixed point cohomogeneity at most one.
\begin{enumerate}[label=(\arabic*),leftmargin=2.2em]
\item Figure~\ref{fig:s6actions}(a). $G=SU(2)$ acts on the first 4 coordinates of $S^6\subset\R^7$ via the realification of $\C^2$ and fixes the last 3 coordinates. The boundary $\partial M^*$ consists of fixed points, and the interior has trivial isotropy. Actions of this type are fixed point homogeneous.
\item $G=SO(3)$, $M^*=B^4$. The group acts on the first 3 coordinates of $S^6\subset\R^7$ via rotation and acts trivially on the last 4 coordinates. The boundary consists of fixed points, and the interior consists of principal orbits with $SO(2)$-isotropy. This is the only case with $\dim(M^*)=4$. Actions of this type are fixed point homogeneous.
\item Figure~\ref{fig:s6actions}(b). $G=SO(3)$ acts on the first 5 coordinates via the irreducible 5-dimensional real representation and acts trivially on the last 2 coordinates. The equator of $\partial M^*$ consists of fixed points, and the two boundary faces corresponding to the two open hemispheres have $O(2)$-isotropy. The interior consists of principal orbits with isotropy $\Z/2\Z\oplus\Z/2\Z$. Actions of this type have fixed point cohomogeneity one.
\end{enumerate}
\end{example}

\begin{example}\label{ex:s6remaining}
Remaining linear actions on $S^6$.
\begin{enumerate}[label=(\arabic*),leftmargin=2.2em]
\item Figure~\ref{fig:s6actions}(c). This action is given by
\[
A(x,y)=(Ax,Ay),\qquad A\in SU(2),\quad x\in\R^4,\quad y\in\R^3,
\]
where the action on $\R^4$ is the realification of the standard representation on $\C^2$ and the action on $\R^3$ is the adjoint representation. Actions of this type are classified by Theorem~\ref{thm:su2main}.
\item Figure~\ref{fig:s6actions}(d). This action is given by
\[
A(x,y,z)=(Ax,Ay,z),\qquad A\in SO(3),\quad x,y\in\R^3,\quad z\in\R.
\]
Actions of this type are classified by Theorem~\ref{thm:so3main}.
\item Let
\[
V=\mathcal H_3(\R^3)=\{p\in\R[x,y,z]\mid p\text{ is homogeneous of degree }3,\ \Delta p=0\}.
\]
This is the irreducible 7-dimensional real representation of $SO(3)$, and the unit sphere $S(V)$ is the round $S^6$. The principal isotropy is trivial. To see this, restrict the representation to rotations around a fixed axis. Its weights are $0,1,2,3$. A rotation of a general angle has a 1-dimensional fixed space, while a half-turn or a rotation through $2\pi/3$ or $4\pi/3$ has a 3-dimensional fixed space. The incidence set of non-identity rotations and their fixed vectors therefore has dimension at most 5 in the 7-dimensional space $V$. Hence a general vector has trivial stabilizer.

The harmonic cubic $p(x,y,z)=xyz$ has stabilizer $A_4$. Indeed, an element preserving its zero set is a signed permutation matrix. If $ge_i=\epsilon_i e_{\sigma(i)}$, then $g\in SO(3)$ gives $\epsilon_1\epsilon_2\epsilon_3\operatorname{sgn}(\sigma)=1$, while $g\cdot xyz=xyz$ gives $\epsilon_1\epsilon_2\epsilon_3=1$. Thus $\sigma$ is even, and the stabilizer has $3\cdot4=12$ elements and is isomorphic to $A_4$. Consequently $SO(3)/A_4$ is an exceptional orbit. The zonal harmonic cubic
\[
z(2z^2-3x^2-3y^2)
\]
has stabilizer $SO(2)$ and gives a singular orbit. This action supplies the tetrahedral example in Theorem~\ref{thm:main}(3).
\end{enumerate}
\end{example}

The effective real representations of $SO(3)$ on $\R^7$ have decompositions
\[
7,\qquad 5+1+1,\qquad 3+3+1,\qquad 3+1+1+1+1,
\]
and the effective $SU(2)$-representations have decompositions $4+3$ and $4+1+1+1$. Thus Examples~\ref{ex:s6low} and \ref{ex:s6remaining} give all effective linear actions under consideration on $S^6$.

\begin{figure}[ht]
\centering
\includegraphics[width=.70\textwidth]{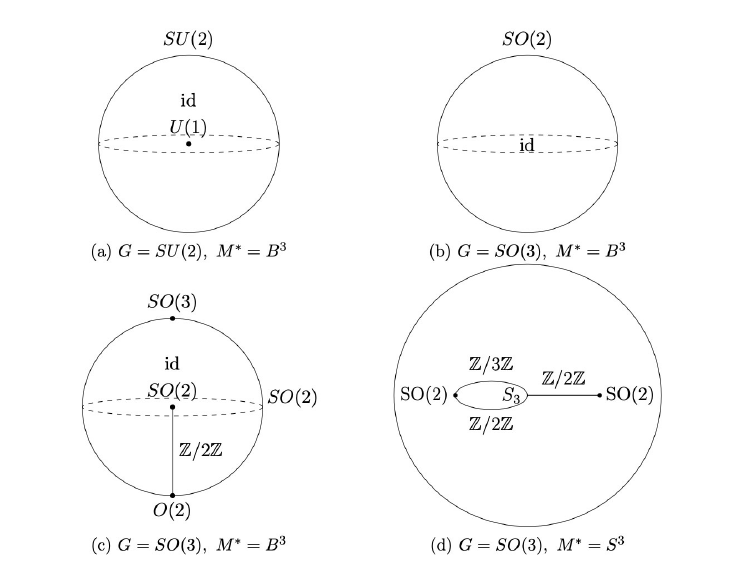}
\caption{$M=\CP^3$.}\label{fig:cp3actions}
\end{figure}

\begin{example}\label{ex:cp3}
Actions on $\CP^3$.
\begin{enumerate}[label=(\arabic*),leftmargin=2.2em]
\item Figure~\ref{fig:cp3actions}(a). A linear $SU(2)$-action on $\CP^3$, acting on the first 2 homogeneous coordinates and fixing the last 2 homogeneous coordinates. We have $M^*=B^3$, $\partial M^*=S^2$ consists of fixed points, and the interior minus the center has trivial isotropy. The center has $U(1)$-isotropy, represented by $[x,y,0,0]\in\CP^3$. Actions of this type are fixed point homogeneous.
\item Figure~\ref{fig:cp3actions}(b). Let $A\in SU(2)$ act on $\CP^3$ via $A(x,y)=(Ax,Ay)$, where $x,y\in\C^2$. This action is ineffective since $-I\in SU(2)$ acts trivially, and hence it descends to an $SO(3)$-action. The interior of $M^*$ consists of principal orbits, and $\partial M^*$ consists of singular $SO(2)$-orbits.
\item Figure~\ref{fig:cp3actions}(c). This action is given by
\[
A(z_1:z_2:z_3:z_4)=(A(z_1,z_2,z_3)^T:z_4),\qquad A\in SO(3).
\]
The orbit space is $B^3$. The interior minus a line segment corresponds to principal orbits, the boundary 2-sphere minus two points corresponds to $SO(2)$-singular orbits, and the two poles correspond to a fixed point and an $O(2)$-orbit. The center of the ball corresponds to an $SO(2)$-orbit, and a line segment in the interior corresponds to $\Z/2\Z$-orbits connecting the $O(2)$-orbit and the interior $SO(2)$-orbit.
\item Figure~\ref{fig:cp3actions}(d). The irreducible representation of $SU(2)$ on $\C^4$ induces an action on $\CP^3$, which has kernel $\Z/2\Z$ and descends to $SO(3)$. For this action, $M^*=S^3$, the principal isotropy is trivial, and there are two singular orbits with $SO(2)$-isotropy. The exceptional orbits are drawn in the picture.
\end{enumerate}
\end{example}

\begin{figure}[ht]
\centering
\includegraphics[width=.68\textwidth]{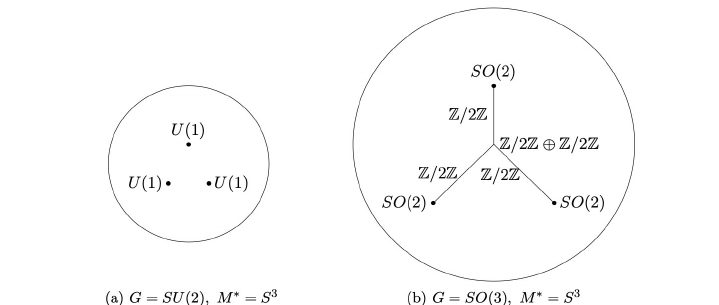}
\caption{$M=SU(3)/T^2$.}\label{fig:flagactions}
\end{figure}

\begin{example}\label{ex:flag}
Actions on $SU(3)/T^2$.
\begin{enumerate}[label=(\arabic*),leftmargin=2.2em]
\item Figure~\ref{fig:flagactions}(a). This action is given by the left multiplication of $SU(2)$ on $SU(3)/T^2$, i.e.,
\[
A\cdot(gT^2)=AgT^2,
\]
where $SU(2)$ is the upper left block in $SU(3)$. The orbit space is a 3-sphere with three singular orbits with $U(1)$-isotropy corresponding to the matrices
\[
I,\qquad
\begin{bmatrix}1&0&0\\0&0&-1\\0&1&0\end{bmatrix},\qquad
\begin{bmatrix}0&1&0\\0&0&1\\1&0&0\end{bmatrix}.
\]
The principal isotropy is trivial. Actions of this type are classified by Theorem~\ref{thm:su2main}.
\item Figure~\ref{fig:flagactions}(b). This action is given by the left multiplication of $SO(3)$ on $SU(3)/T^2$. We have $M^*=S^3$, and there are three singular orbits. The orbit strata are indicated in the picture.
\end{enumerate}
\end{example}

\begin{figure}[H]
\centering
\includegraphics[width=.36\textwidth]{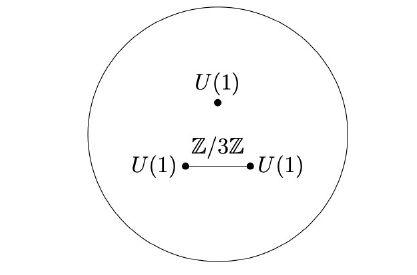}
\caption{$G=SU(2)$, $M^*=S^3$ on $SU(3)//T^2$.}\label{fig:biquotient}
\end{figure}

\begin{example}\label{ex:biquotient}
An action on $SU(3)//T^2$. Recall that
\[
T^2=\{(\diag(z,w,zw),\diag(1,1,\bar z^2\bar w^2))\mid z,w\in S^1\}
\subset S(U(3)\times U(3))
\]
acts freely on $SU(3)$. The group $G=SU(2)$ acts from the right as the upper left block of $SU(3)$, commuting with the $T^2$-action. We have $M^*=S^3$. A computation using the Mayer--Vietoris sequence shows that $G$-spaces of this type have the same cohomology groups as $SU(3)//T^2$, that is,
\[
H^0=H^6=\Z,\qquad H^2=H^4=\Z\oplus\Z,\qquad H^{2i+1}=0.
\]
The orbit stratification is depicted in Figure~\ref{fig:biquotient}.
\end{example}

\Needspace{10\baselineskip}
\medskip
\noindent\textbf{Conflict of Interest}\quad The author declares no conflict of interest.

\medskip
\noindent\textbf{Acknowledgements}\quad This paper is part of the author's Ph.D. thesis under the supervision of Professor Wolfgang Ziller and we would like to thank him for his endless advice and support. We would also like to thank Fuquan Fang, Francisco Gozzi, Karsten Grove, Xiaochun Rong and Fabio Simas for helpful conversations. Lastly we thank the referees for their time and valuable comments.

\end{document}